\documentclass[12pt]{amsart}
\usepackage[dvips]{graphicx}
\newcommand{\C}{{\mathbb{C}}}

\renewcommand{\H}{{\mathbb{H}}}
\newcommand{\Z}{{\mathbb{Z}}}
\newcommand{\N}{{\mathbb{N}}}
\newcommand{\Q}{{\mathbb{Q}}}
\newcommand{\R}{{\mathbb{R}}}

\newlength{\figboxwidth}
\def\defeq{\ensuremath{\stackrel{\mathrm{def}}{=}}}
\newcommand{\subgroup}{<}

\newcommand{\normal}{\triangleleft}

\newcommand{\arrow}{\rightarrow}

\newcommand{\compos}{\circ}

\newcommand{\isom}{\cong}

\newcommand{\superset}{\supset}
\newcommand{\tensor}{\otimes}
\newcommand{\bs}{\backslash}

\newcommand{\Aut}{\operatorname{Aut}}

\newcommand{\girth}{\operatorname{girth}}

\renewcommand{\mod}{\operatorname{mod}}

\newcommand{\Stab}{\operatorname{Stab}}

\newcommand{\Lk}{\operatorname{Lk}}

\newtheorem{theorem}{Theorem}[section]

\newtheorem{example}[theorem]{Example}
\newtheorem{proposition}[theorem]{Proposition}
\newtheorem{lemma}[theorem]{Lemma}

\newtheorem{definition}[theorem]{Definition}
\newtheorem{question}[theorem]{Question}

\mathchardef\GG="321D
\input xy
\xyoption{all}
\numberwithin{equation}{section}

\def\defeq{\ensuremath{\stackrel{\mathrm{def}}{=}}}

\title{Ramanujan graphs with small girth.}
\author{Yair Glasner}
\email{yair@math.uic.edu}   
\urladdr{http://www.math.uic.edu/~yair} 
\address{Department of Mathematics, University of Illinois at Chicago, 
        851 South Morgan Street, Chicago, IL, 60607-7045, USA}
\date{}
\begin{document}
\bibliographystyle{amsalpha}
\begin{abstract}
    We construct an infinite family of $(q+1)-$regular Ramanujan graphs $X_n$ of
    girth $1$. We also give covering maps $X_{n+1} \arrow X_n$
    such that the minimal common covering of all the $X_n$'s is the
    universal covering tree. 
\end{abstract}
\maketitle

%-----------------------------------------------------------------------------

\section{Introduction}
    Ramanujan graphs were first introduced by Lubotzky Phillips and Sarnak (LPS) in
    \cite{LPS:ramanujan} as graphs satisfying the ``asymptotically optimal''
    bound on the size of the second eigenvalue. A $(q+1)-$regular
    graph $X$ is called {\it{Ramanujan}} if $\left| \lambda \right| \le 2\sqrt{q}$
    for every non-trivial (i.e. $\ne \pm(q+1)$) eigenvalue $\lambda$
    of its adjacency matrix. By saying that the bound $2\sqrt{q}$ is
    {\it{asymptotically optimal}} we mean that trying to impose
    any lower bound gives rise only to finite families of graphs; this is the
    content of Alon-Boppana theorem \cite{Alon:AB}. It is an open question weather there
    exist infinite families of $(q+1)-$regular Ramanujan graphs for a general
    number $q$. Explicit constructions in the case where $q$ is prime were
    given by Lubotzky Phillips and Sarnak (LPS graphs),
    these were later extended to the case where
    $q$ is a prime power (\cite{Morgenstern:Ramanujan,Jordan_Livne:graphs})
    \footnote{It should be emphasized that even non constructive methods or methods of
          probabilistic nature for proving existence of Ramanujan graphs are not
          known. The best known results so far are due to Joel Friedman
          (\cite{Friedman:second_ev}) who shows that the non-trivial eigenvalues
          of a random $(q+1)-$regular graph with $n$ vertices satisfy
          $\left| \lambda \right| \le 2\sqrt{q} + \log(q) + C$ almost surely
              (i.e. with a probability that tends to 1 with $n$). The bound here
              is strictly larger than the Ramanujan bound and independent of $n$.}.
    It is not difficult to describe these examples, but the proof of the Ramanujan
    property relies on very deep theorems from number theory.

    LPS graphs exhibit many interesting combinatorial properties, some of these
    are a direct consequence of the Ramanujan property and others are independent of
    the spectral properties of the graph. Two examples are:
    \begin{itemize}
      \item It is a direct consequence of the Ramanujan
            property that LPS graphs are good expanders.
      \item It can be proved in an elementary way, independent of the Ramanujan property,
        that LPS graphs have very large girth. In fact the bi-partite LPS graphs
        satisfy $\girth(X) \ge \tfrac{4}{3}\log(|X|)$.
    \end{itemize}
    Lubotzky, in his book \cite[Question 10.7.1]{Lubotzky:Book}, poses
    the question of clarifying the connection between the Ramanujan property
    and the girth. There are some theorems showing a correlation between
    the eigenvalue distribution and the existence of small circuits, but they are all
    rather weak.
    For example Greenberg (see \cite{Greenberg},\cite[theorem 4.5.7]{Lubotzky:Book})
    proves that an infinite family of Ramanujan graphs
    $\ldots X_n \arrow X_{n-1} \arrow \ldots \arrow X_1$ covering each other,
    with $X_n \arrow X_1$ a regular covering map,
    must satisfy $\girth(X_n) \arrow \infty$
    (In other words an example like the one given in theorem (\ref{thm:main}) below is
    not possible when the covering maps $X_n \arrow X_1$ are assumed to be
    regular).
    In the other direction McKay
    (see \cite{McKay:eigenvalues})
    shows that infinite families of $(q+1)-$regular graphs with asymptotically few circuits
    have most of their eigenvalues concentrated in the Ramanujan interval
    $[-2\sqrt{q},2\sqrt{q}]$.

    In this paper we give the following example, proving that the Ramanujan
    property does not imply large girth.
    \begin{definition}
      An infinite sequence of graphs covering each other
      $\ldots \arrow X_{n+1} \arrow X_n \arrow X_{n-1} \arrow \ldots \arrow X_1$ will
      be called a {\it{tower of graphs}}.
    \end{definition}

    \begin{theorem}
    \label{thm:main}
    There exists a tower of regular Ramanujan graphs $X_n$
    satisfying the following properties:
    \begin{enumerate}
       \item \label{itm:girth}
        There is a common bound $M$ on the girth of all the graphs.
       \item \label{itm:inv}
        The minimal common covering of all graphs is the universal covering tree.
    \end{enumerate}
    \end{theorem}
    \noindent

    We also give an explicit description of one such family of graphs, similar
    to the explicit constructions given by Lubotzky Phillips and
    Sarnak in (\cite{LPS:ramanujan}).
    In fact in section (\ref{sec:concrete}) we show that our graphs can be realized
    as families of Schreier graphs of $PSL_2(\Z/q^n\Z)$ (or closely related
    groups) with respect to some Cartan subgroups. The Ramanujan graphs obtained in
    these examples have girth $1$, (i.e. they all contain loops).

    There are good reasons to expect that examples such as the one given in theorem
    \ref{thm:main} should exist. Indeed, the property of being a Ramanujan graph
    is only asymptotically optimal. It is not difficult to find small graphs,
    satisfying better bounds on the second eigenvalue. Given such
    a ``better than Ramanujan'' graph we can make a local change to the graph:
    creating a small circuit while introducing only minor changes in the eigenvalues
    and retaining the Ramanujan property. To obtain an infinite family
    of Ramanujan graphs with small girth, we have to start with an infinite family
    of better then Ramanujan graphs with some precise estimates on their second eigenvalues,
    which seems very difficult.
    We do not take this approach, instead
    we go back and introduce a minor change in the construction of LPS graphs and
    this yields the desired family of Ramanujan graphs with small girth.

    If theorem \ref{thm:main}(\ref{itm:inv}) does not hold then
    \ref{thm:main}(\ref{itm:girth}) automatically does hold because any closed path
    in the common covering of all Ramanujan graphs will appear in each one of them.
    This happens exactly when the fundamental groups of our graphs have a non-trivial
    intersection. In \cite{LPS:ramanujan} Ramanujan graphs are constructed
    as quotients of the Bruhat-Tits tree $T$ of $PGL_2(\Q_q)$ by torsion free
    congruence subgroups of a $\{q\}$-arithmetic lattice $\Gamma \subgroup PGL_2(\Q_q)$.
    It is customary to use principal congruence subgroups but, as the intersection of
    an infinite family of such is always trivial, we replace them by an infinite
    family of non-principal congruence subgroups. In order to achieve a minimal common
    covering which is a tree (\ref{thm:main}(\ref{itm:inv})), we change the covering
    morphisms $X_{n+1} \arrow X_n$.

    {\bf{Remark:}}
      We start from the classical construction of Ramanujan
      graphs due to Lubotzky Phillips and Sarnak and modify it slightly.
      This construction involves
      definite quaternion algebras defined over $\Q$, and is described in
      \cite{Lubotzky:Book} and in section (\ref{sec:LPS}).
      This approach has the advantage of simplifying the presentation and the disadvantage
      of yielding only $(q+1)-$regular graphs, where $q$ is prime.
      The regularity restriction is not really necessary, our
      method will work also for variants of the LPS construction
      which yield $(q+1)-$regular Ramanujan graphs for every prime power $q$.
      For example Morgenstern's construction in positive
      characteristic \cite{Morgenstern:Ramanujan} or the construction of Jordan Livn{\'e}
      involving totally definite quaternion algebras over number fields
      \cite{Jordan_Livne:graphs}.

    I would like to thank Shlomo Hoory, Natan Linial, Ron Livn{\'e}, Alex Lubotzky,
    Shahar Mozes, and  Andrezej {\.Z}uk for many helpful discussions on Ramanujan graphs.
    I also thank Alex Lubotzky for showing me the proof of lemma
    \ref{lem:AbelianFG} and Shahar Mozes for letting me include his proof of Prasad's theorem
    (\ref{prop:torus}). Finally, I wish to thank the two referees. They read carefully the
    original manuscript, and their many comments helped me bring it to its current
    more readable form.

\section{Review}
\label{sec:review}
\subsection{Covering theory for graphs}
    We define a {\it{graph}} $X$ to be a set of {\it{vertices}} $VX$, a set of {\it{edges}}
    $EX$ together with a fixed-point-free involution $^{-}:EX \arrow EX$
    (associating with every edge $e$ an edge $\overline{e}$ called its inverse)
    and two maps called {\it{the origin and terminus maps}} $o,t:EX \arrow VX$ satisfying
    $o\overline{e} = te, \quad t\overline{e} = oe$. We think of our graphs as
    non-directed but we represent each {\it{geometric edge}} $[e]$
    by a pair of directed edges $\{e,\overline{e}\}$. All our graphs
    might contain loops or multiple edges. All the graphs in this paper are assumed
    connected. Notions like graph {\it{morphisms, automorphisms, etc \ldots}} are
    all defined in the obvious way. All automorphisms of a given graph $X$ form
    a group denoted by $\Aut(X)$.

    Let $G$ be a group acting (on the left, by graph automorphisms) on a graph $X$.
    We say that the action is {\it{without inversion}} if no element of $G$
    takes an edge to its inverse. Whenever $G$ acts without inversion on a graph
    $X$ there exists a well defined {\it{quotient graph}} denoted $G \bs X$ and a
    {\it{quotient morphism}} $p:X \arrow G \bs X$.
    To see this one should check that the structure maps $\{^{-},o,t\}$ of the graph $X$
    induce well defined structure maps on the quotient sets $G \bs VX$ and $G \bs EX$.
    We say that a group $G$ {\it{acts freely}} on the graph $X$, if all vertex stabilizers
    are trivial.

    A {\it{link of a vertex}}
    $v \in VX$ is the set $\Lk(v) \defeq \{ e \in EX | oe =v \}$.
    A graph morphism $\phi:X \arrow Y$ induces a map
    $\phi_v:\Lk(v) \arrow \Lk(\phi(v)) \ \forall v \in VX $.
    A {\it{covering map of graphs}} $\phi:X \arrow Y$ is, by definition, a (surjective)
    map that induces a bijection on every vertex link.
    If a group $G$ acts freely and without inversion on a graph $X$,
    the quotient map $p: X \arrow G \bs X$ is a covering
    map. Covering maps obtained in this way are called {\it{regular}}.
    Every graph is regularly covered by a tree:
    \begin{theorem} \label{thm:universal_covering}
       For every graph $Y$ there exists a unique regular covering map
       $\pi:\tilde{Y} \arrow Y$ with $\tilde{Y}$ a tree.
    \end{theorem}
    The tree $\tilde{Y}$ is called
    {\it{the universal covering tree of $Y$}}. The group that acts freely and without
    inversion to yield $Y$ as a quotient is denoted by $\pi_1(Y,\cdot)$, it is
    called {\it{the fundamental group of $Y$}}. The uniqueness statement in the
    theorem means that
    the pair $(\pi_1(Y,\cdot),\tilde{Y})$ is uniquely determined by the graph $Y$
    (up to a naturally defined notion of isomorphism of group actions on graphs).
    It turns out that fundamental groups of graphs are always free groups.

    If a group $G$ acts freely and without inversion on a graph $X$ with a quotient
    graph $Y=G \bs X$, and if $H<G$
    is a subgroup there is a natural covering map of graphs
    $X \arrow H \bs X \arrow G \bs X = Y$. It turns out that if we take $X = \tilde{Y}$
    and $G = \pi_1(Y,\cdot)$ all possible covering graphs of $Y$ are obtained in this way:
    \begin{theorem} \label{thm:Galois}
       There is a bijective correspondence, called {\it{the Galois correspondence}},
       between the subgroups of $\pi_1(Y,\cdot)$ and intermediate covering
       graphs $\tilde{Y} \arrow Z \arrow Y$. The Galois correspondence associates
       with a subgroup $H$ of $\pi_1(Y,\cdot)$ the intermediate covering:
       \begin{equation} \label{eqn:Galois}
        \tilde{Y} \arrow H \bs \tilde{Y} \arrow
        \left(Y = \pi_1(Y,\cdot) \bs \tilde{Y}\right).
       \end{equation}
       Inclusion of subgroups $H_2 \subgroup H_1 \subgroup \pi_1(Y,\cdot)$ is transformed
       to covering of graphs $H_2 \bs \tilde{Y} \arrow H_1 \bs \tilde{Y}$.
       Normal subgroups of $\pi_1(Y,\cdot)$ correspond to regular coverings of $Y$.
       Explicitly: if $N \normal \pi_1(Y,\cdot)$
       the action of $\pi_1(Y,\cdot)$ on $\tilde{Y}$ induces a free action
       without inversion of $\pi_1(Y,\cdot) / N$ on $N \bs \tilde{Y}$; and
       there is a natural isomorphism
       $(\pi_1(X,\cdot) / N) \bs (N \bs \tilde{X}) \cong Y$.
    \end{theorem}

    The $(q+1)-$regular graphs are exactly the graphs whose universal
    covering tree is $T=T_{q+1}$, the $(q+1)-$regular tree. Let $q+1 = 2n$ and let
    $X$ be the graph with one vertex and $n$ edges, such a graph
    is sometimes called a {\it{wedge of $n$ circles}}.
    The Galois correspondence can be described very explicitly for the graph $X$:
    \begin{example} \label{ex:rose} \
    \begin{itemize}
    \item \label{itm:general}
    Let $\Gamma$ be a free group on $n$ letters,
    $S = \{\gamma_1,\gamma_2,\ldots,\gamma_n,
         \gamma_1^{-1},\gamma_2^{-1},\ldots,\gamma_n^{-1}\} \subset \Gamma$ a set of
    free generators and their inverses,
    and $T = X(\Gamma,S)$ the right Cayley graph of $\Gamma$ with respect to $S$.
    $T = T_{2n}$ is the $2n-$regular tree and $\Gamma$ acts freely transitively and without
    inversion on $T$, thus $X = \Gamma \bs T$ is a wedge of $n$ circles and the pair
    $(\Gamma,T)$ can be identified with $\left( \pi_1(X,\cdot),\tilde{X} \right)$.
    \item \label{itm:Cayley}
    If $N \normal \Gamma$, the graph
    $N \bs T$ can be identified with the (right) Cayley graph of $\Gamma/N$
    with respect to the symmetric set of generators
    $\{\overline{\gamma_1},\overline{\gamma_2},\ldots,
    \overline{\gamma_n}^{-1} \}$
    \item \label{itm:Schreier}
    If $N \subgroup H \subgroup \Gamma$
    are subgroups with $N \normal \Gamma$, the graph
    $H \bs T$ can be identified with the Schreier graph of the group pair
    $(\Gamma/N , H/N)$ with respect to the same set of generators.
    Explicitly: we can identify the vertices of $H \bs T$ with the right
    cosets of $(H/N)$ in $(\Gamma/N)$. The directed edges will
    be of the form
    $\{(\gamma N,\gamma \gamma_i^{\pm 1} N)\}_{\gamma \in \Gamma, \gamma_i^{\pm 1} \in S}$.
    \end{itemize}
    \end{example}

\subsection{Classification of tree automorphisms}
    \begin{theorem} \label{thm:tree_aut} (See \cite{Serre:Trees})
    Let $T$ be a regular tree, $d$ the standard metric on $VT$,
    $\sigma \in Aut(T)$,
    $l = l(\sigma) = \min\{d(v,\sigma v):v \in VT\}$ and
    $X=X(\sigma) = \{v \in VT|d(v,\sigma v) = l(\sigma)\}$,
    then exactly one of the following three possibilities hold:
    \begin{enumerate}
    \item \label{itm:inversion}
        {\underline{$\sigma$ is an {\it{inversion}}:}}
        There exists an edge $e \in ET$ with
        $\sigma e = \overline{e}$. in this case $e$ is uniquely determined,
        $l = 1$ and $X$ consists of the two points $\{oe,te\}$.
    \item \label{itm:elliptic}
        {\underline{$\sigma$ is {\it{elliptic}}:}}
        $l = 0$ and $X$ is a convex subset of $T$ consisting of all the
        points fixed by $\sigma$.
    \item \label{itm:hyperbolic}
        {\underline{$\sigma$ is {\it{hyperbolic}}:}} $l \ge 0$ and $X$
        is a bi-infinite line on which $\sigma$ acts by a translation of length $l$.
    \end{enumerate}
    \end{theorem}

\subsection{Bruhat-Tits theory} \label{sec:BT}
    A rich source of examples for group actions on trees comes from Bruhat-Tits
    theory. The basic tool is a natural action of the group
    $G=PGL_2(\Q_q)$ on the $(q+1)-$regular tree $T = T_{q+1}$
    (an excellent exposition of this action can be found in Serre's book \cite{Serre:Trees}).
    We will use the following properties of this action:
    \begin{enumerate}
    \item \label{itm:stabilizers}
        The subgroup $K = PGL_2(\Z_{q})$ is the stabilizer of a vertex
        $O \in VT$.
        $G$ acts transitively on the vertices and on the directed edges of $T$.
        In particular all vertex stabilizers are conjugate to $K$.
    \item \label{itm:spheres}
        The stabilizer of the sphere $S_{T}(O,n)$ is the group
        $K(q^{n}) \defeq ker\{PGL_2(\Z_{q}) \stackrel{\psi_n}{\arrow}
                 PGL_2(\Z/q^n\Z)\}$.
        Thus the action of $K$ on $S_{T}(O,n)$ factors through an action
        of $PGL_2(\Z/q^n\Z)$ on $S_{T}(O,n)$.
        It turns out that this action can be identified with the
        action of $PGL_2(\Z/q^n\Z)$ on the projective line ${\mathbb{P}}^{1}(\Z/q^n\Z)$
        by M\"obius transformation.
        Consequently, the action of $PGL_2(\Z/q^n\Z)$ on $S_{T}(O,n)$
        is transitive and the point stabilizers are all conjugate to the
        (Borel) subgroup of upper triangular matrices. Since the action of
        $PGL_2(\Z/q^n\Z)$ on $S_{T}(O,n)$ comes from an action on the ball
        $B_{T}(O,n)$ it is impossible for this action to be $2-$transitive,
        it is transitive however on the pairs of points
        $\{(x,y) \in S_{T}(O,n) \times S_{T}(O,n)| \text{ The path [x,y] contains
        the point } O\}$ and the stabilizers of such pairs are conjugate to the
        diagonal (Cartan) subgroup \footnote{In fact the action is also transitive
        on the set of triplets of points satisfying a similar geometric condition, the
        stabilizer of such a triplet is trivial}.
    \item \label{itm:boundary}
        The set of infinite rays emerging from $O$, is called the {\it{boundary of the
        tree}} and the group $K$ acts transitively on it.
        In fact this action on the boundary can be identified with the action
        of $PGL_2(\Z_q)$ on the projective line
        ${\mathbb{P}}^{1}(\Z_q) \cong {\mathbb{P}}^{1}(\Q_q)$ by M\"obius
        transformations
        \footnote {This is proved by noticing that everything that was said about
                   the action on the sphere $S(n)$ is compatible with the natural
               map $S(n+1) \arrow S(n)$ and then passing to the inverse limit.}.
    \end{enumerate}

\subsection{The construction of LPS graphs}
\label{sec:LPS}
    By theorem (\ref{thm:Galois}) $(q+1)-$regular graphs are equivalent
    to groups $\Gamma$ acting freely and with finitely many orbits on the
    $(q+1)-$regular tree $T_{q+1}$.
    A wider class of groups is the class of {\it{uniform lattices}}
    namely groups acting with a finite number of orbits and with finite
    vertex stabilizers. One method for constructing uniform lattices inside
    $PGL_2(\Q_q)$
    is the arithmetic construction due to Borel and Harish-Chandra
    (see \cite{BHC:Arithmetic_I,BHC:Arithmetic_II}).
    Uniform lattices constructed in this way are referred to as
    $\{q\}-$arithmetic lattices or sometimes just as arithmetic lattices.
    The extensive knowledge available on arithmetic lattices and their properties has
    enabled Lubotzky Phillips and Sarnak
    (\cite{LPS:ramanujan}) to prove that $\Gamma(n) \bs T$ are Ramanujan
    graphs for certain infinite families of arithmetic lattices $\{\Gamma(n)\}_{n \in \N}$.
    We proceed to describe these lattices.

    Let $\H = \H_{u,v}$ be a {\it{quaternion algebra}}
    defined over $\Q$ (i.e. $u,v \in \Q$).
    $\H$ is a four dimensional algebra over $\Q$ spanned by the four
    symbols $1,i,j,k$ satisfying the relations $i^2= -u, j^2= -v; \ ij = -ji =k$.
    If $k/\Q$ is a field, $\H(k) \defeq \H \tensor_{\Q} k$ (or in other words $\H(k)$
    is defined in the same way only as an algebra over $k$).
    Let $G$ be the $\Q$-algebraic group $\H^{*}/Z\H^*$, one can think of $G$
    as a way of associating a group $G(k)$ to every field $k/\Q$ by letting
    $G(k)=\H(k)^{*}/Z\H(k)^*$ (The group of invertible elements modulo the center).
    If $R < \C$ is a ring we will think of
    $G(R)$ as the group of elements in $G(\C)$ that have representatives all of whose
    coefficients are elements of $R$.

    There is a dichotomy saying that $H(k)$ is either isomorphic
    to $M_2(k)$ or is a division algebra
    \footnote{It is easy to see that $\H$ is a central simple algebra and from the
    structure theory of such algebras $\H(k) \cong M_n(D)$ for some division ring
    $D/k$, the dichotomy now follows because $dim_{k} \H(k) = 4$.}.
    In the first case we say that $\H$ {\it{splits}} over $k$ and in the
    second that $\H$ {\it{ramifies}} over $k$. If $\H$ splits over $k$ then
    $G(k) \cong PGL_2(k)$.

    We say that $\H$ splits (resp. ramifies)
    at the prime $q$ if it splits (resp. ramifies) over $\Q_q$, (including the case
    $q = \infty, \ \Q_q = \R$). Every quaternion algebra $\H$ ramifies over a
    finite set of primes (of even cardinality) and splits over all other primes.
    $\H$ is called a {\it{definite}} quaternion algebra if it ramifies at $\infty$.

    Assume that a definite quaternion algebra $\H$ splits at $q$, then
    $G(\Z[1/q]) < G(\Q_{q}) \cong PGL_2(\Q_{q})$ is an example of a $\{q\}-$
    arithmetic lattice in $PGL_2(\Q_{q})$.
    If $N$ is any integer prime to $q$ the homomorphism
    $\psi_N:\Z[1/q] \arrow \Z[1/q]/N\Z[1/q] \isom \Z/N\Z$ gives rise to a
    homomorphism
    \begin{equation}
        \psi_N: \Gamma = G(\Z[1/q]) \arrow G(\Z[1/q]/N\Z[1/q]) \isom PGL_2(\Z/N\Z).
    \end{equation}
    \begin{definition}
        Principal congruence subgroups of $\Gamma$ are groups of the form
        $\Gamma(N) = ker(\psi_N)$. A subgroup of $\Gamma$ is called a
        congruence subgroup if it contains
        a principal congruence subgroup.
    \end{definition}
    \noindent
    If $N$ is big enough then the group $\Gamma(N)<PGL_2(\Q_{q})$
    acts freely on the
    Bruhat-Tits tree $T_{q+1}$.
    \begin{theorem} \label{thm:LPS} (Lubotzky, Phillips, Sarnak).
        Let $\H = \H_{u,v}$ be a definite quaternion algebra defined over $\Q$,
        $G$ the $\Q-$algebraic group $\H^{*}/Z\H^{*}$ and assume that $\H$ splits
        over $q$. If a congruence subgroup $\Gamma' < \Gamma = G(\Z[1/q])$
        acts freely on $T = T_{q+1}$. then
        the graph $\Gamma' \bs T$ is a $(q+1)-$regular Ramanujan graph.
    \end{theorem}
    {\bf{Remark: }}
            Theorem (\ref{thm:LPS}) is usually stated for principal
        congruence subgroups. The general case is easily deduced:
        by definition every congruence subgroup $\Gamma'$
        contains a principal congruence subgroup $\Gamma(N)$ giving rise to
        a covering map $\eta:\Gamma(N) \bs T \arrow \Gamma' \bs T$. If $f$ is an
        eigenfunction for the adjacency operator on $\Gamma' \bs T$ then
        $f \compos \eta$ will be an eigenfunction on $\Gamma(N) \bs T$ with the
        same eigenvalue. The graph
        $\Gamma(N) \bs T$ will, therefore, inherit all the eigenvalues of
        $\Gamma' \bs T$ proving that the later must be Ramanujan if the former is.

\section{Ramanujan Graphs of Small Girth}
    Let $\H = \H_{u,v}$ be a definite quaternion algebra defined over $\Q$
    ($0 < u,v \in \Q$, $i^2= -u,j^2=-v,ij=-ji=k$, $\H$ ramifies at $\infty$).
    Let $G$ be the $\Q$ algebraic group $\H^{*}/Z\H^{*}$, $q_1$ a prime such
    that $\H$ splits at $q_1$ and $\Gamma$ the $\{q_1\}-$arithmetic lattice
    $G(\Z[1/q_1]) < PGL_2(\Q_{q_1})$.

    We give a geometric description of some congruence subgroups of $\Gamma$
    which is based on ideas coming from the theory of lattices acting on products of
    trees (\cite{BM:Lattices,Jordan_Livne:complexes}).
    Let $q_2 \ne q_1$ be a second prime such that $\H$ splits over $q_2$
    and $T_{2} = T_{q_2+1}$ the Bruhat-Tits tree corresponding to $PGL_2(\Q_{q_2})$.
    The group $\Gamma$ acts on $T_2$ through its embedding in
    $G(\Q_{q_2}) \cong PGL_2(\Q_{q_2})$, in fact here $\Gamma$ is a subgroup of
    $PGL_2(\Z_{q_2})$ so it fixes a vertex $O_2 \in VT_2$.
    Furthermore by \ref{sec:BT}(\ref{itm:spheres}) we can identify
    $\Gamma(q_2^n) = \Gamma_{B_{T_2}(O_2,n)}$.

    Given any finite set of vertices $O_2 \in C \subset VT_2$ there exist an $n \in \N$
    such that $\Gamma_C \superset \Gamma_{B_{T_2}(O_2,n)}$.
    By \ref{sec:BT}(\ref{itm:spheres}) $\Gamma_{B_{T_2}(O_2,n)} = \Gamma(q_2^n)$ is
    a principal congruence subgroup,
    so the pointwise stabilizer $\Gamma_C < \Gamma$ is a congruence subgroup
    of $\Gamma$, and the corresponding graph $X_C \defeq \Gamma_C \bs T_1$ is
    Ramanujan by theorem (\ref{thm:LPS}). If we have two such sets $C_1 \subset C_2$
    then $\Gamma_{C_2} < \Gamma_{C_1}$ and by the Galois correspondence \ref{thm:Galois}
    there is a covering map $X_{C_2} \arrow X_{C_1}$.

    An ascending sequence of finite subsets
    $O_2 \in C_1 \subset \ldots C_{n-1} \subset C_n \subset C_{n+1} \ldots$
    gives rise a tower of Ramanujan graphs $\{X_{n} \defeq X_{C_n}\}_{n \in \N}$.
    If $C = \bigcup_n C_n \subset VT_2$, then
    $\bigcap_n \Gamma_{C_n} = \Gamma_C$.
    The minimal graph that covers all the $X_n$'s is $X_C \defeq \Gamma_C \bs T_1$. it
    will be a tree and
    theorem \ref{thm:main}(\ref{itm:inv}) will be satisfied iff
    $\Gamma_C = \langle e \rangle$.

    We first construct a tower of Ramanujan graphs whose minimal common covering is not
    a tree and consequently the girth of all graphs is bounded.
    Let $\Sigma$ be the $\{q_1,q_2\}-$arithmetic lattice
    $\Sigma \defeq G(\Z[1/q_1,1/q_2])$ and
    identify $\Gamma$ with the subgroup of $\Sigma$ consisting of these elements
    whose coefficients contain only denominators which are powers of $q_1$, which is
    exactly the subgroup of $\Sigma$ fixing the vertex $O_2 \in VT_2$.
    \begin{equation}
        \Gamma = \Sigma \cap PGL_2(\Z_{q_2}) = \Sigma_{O_2}
    \end{equation}
    Without loss of generality, we may assume that $\Gamma$ acts freely on $T_1$ by
    replacing both $\Sigma$ and $\Gamma$ by $\Sigma(N)$ and $\Gamma(N)$ - the principal
    congruence subgroups mod $N$, where $N$ is a large enough number such that
    $(N,q_1q_2)= 1$.

    \begin{definition}
       We say that the action of the group $\Sigma$ on $T_1 \times T_2$ contains
       a {\it{torus}} if there are two infinite geodesics (=infinite
       paths without backtracking) $l_i \subset T_i$
       and a subgroup $\Z^2 \cong \langle \gamma,\delta \rangle \subgroup \Sigma$
       fixing setwise the tessellated plain $A = l_1 \times l_2 \subset T_1 \times T_2$,
       and acting on it freely and co-compactly.
        \end{definition}
    The action of the group $\Sigma$ on $T_1 \times T_2$ does contain a torus.
    This is a theorem due to Prasad (see \cite{Prasad:Lattices_local_fields}),
    another, more geometrical, proof for the existence of a torus is given
    by Mozes in (\cite[see the discussion below prposition 2.11]{moz:cartan}).
    Both proofs are given in a much greater generality. Shahar Mozes indicated
    to me that his proof, when it is adopted to our particular case, becomes very simple.
    I sketch his argument here for the convenience of the readers.
    In the specific example, described in section (\ref{sec:concrete}), it is easy
    to find a torus explicitly.
    \begin{proposition} (Prasad)
    \label{prop:torus}
    The action of the group $\Sigma$ on $T_1 \times T_2$ contains a torus.
    \end{proposition}
    \begin{proof}
        The reader should refer to figure (\ref{fig:torus}) to understand this
        proof.
        \begin{figure}
        \begin{center}
                \includegraphics{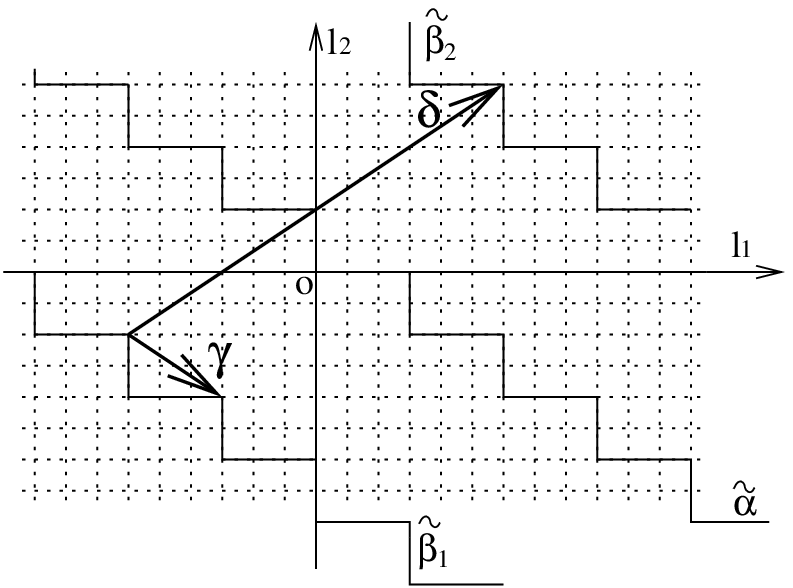}
        \end{center}
        \caption{A torus} \label{fig:torus}
        \end{figure}
        Consider the square complex $X \defeq \Sigma \bs (T_1 \times T_2)$,
        let $p: T_1 \times T_2 \arrow X$ be the covering morphism and
        $O = (O_1,O_2) \in T_1 \times T_2$ a base vertex.
        Choose two closed paths $\alpha_1$ and $\alpha_2$ in the horizontal and
        vertical $1-$skeletons of the complex $X$ respectively. Both $\alpha_i$ should
        start at the base vertex $p(O) \in X$. Now lift the bi-infinite path
        $\alpha \defeq \ldots \alpha_1 \cdot \alpha_2 \cdot \alpha_1 \cdot \alpha_2 \ldots$
        to a path $\tilde{\alpha}$ in $T_1 \times T_2$ passing through $O$.
        The convex hull of the path $\tilde{\alpha}$
        in $T_1 \times T_2$ is a tessellated plane which we denote by $A$,
        it is obviously invariant under the deck transformation
        $\gamma \defeq \alpha_1 \cdot \alpha_2$,
        viewed as an element of $\pi_1(X,p(O))$ acting on $T_1 \times T_2$.
        Now consider all the zig-zag lines parallel to $\tilde{\alpha}$. By the
        invariance under $\gamma$, the
        restriction of $p$ to such a line is determined by a finite segment of a
        fixed length.
        So it is possible to find two such lines $\tilde{\beta_1},\tilde{\beta_2}$
        such that $p|_{\tilde{\beta_1}} = p|_{\tilde{\beta_2}}$. Let
        $\delta : A \arrow A$ be an affine transformation taking
        $\tilde{\beta_1}$ to $\tilde{\beta_2}$. Since $A$ is the convex hull
        of each of the $\tilde{\beta_i}$, the mapping $p|_{A}$ is determined
        by its restriction to $\beta_{i}$. This implies that $\delta$ is also a
        deck transformation $p \compos \delta = p$. The desired $\Z^{2}$ is now generated
        by $\Z^{2} = \langle \gamma,\delta \rangle$ and we are done.
\end{proof}

    We may assume without loss of generality that the action
    $\Z^{2} \cong \langle \gamma,\delta \rangle \leadsto A$
    is the standard action: $\gamma$ acts by translation of length
    $M_1 = l^{T_1}(\gamma)$ on the first coordinate and fixes the second coordinate
    and the other way around for $\delta$.
    $\gamma$ will thus be elliptic on $T_2$ and hyperbolic
    on $T_{1}$ and the opposite for $\delta$. We can also assume that
    $\gamma \in \Gamma$ i.e. that $O_2$, the vertex fixed by $PGL_2(\Z_{q_2})$,
    is contained in $l_2$.

    The element $\gamma$,
    acting on $T_2$, fixes pointwise the bi-infinite line
    $l_2 \in T_2$, in particular it fixes pointwise any finite subset of $l_2$.
    Taking $\{O_2 \in C_n \subset l_2\}_{n \in \N}$ to be symmetric line
    segments of length $2n$ around $O_2$ we obtain a tower of Ramanujan graphs
    $X_n \defeq X_{C_n} = \Gamma_{C_n} \bs T_1$ with bounded girth.
    In fact we can explicitly describe the circle common to all these graphs,
    the infinite line $l_1 \subset T_1$ is mapped to a circle of length $\le M_1$
    in the graph $X_C = \Gamma_{C} \bs T_1$ which in turn covers all the finite graphs $X_n$.

    We now modify our example so that \ref{thm:main}(\ref{itm:inv}) is satisfied,
    i.e so that the minimal common covering of the tower of
    Ramanujan graphs is the tree $T_1$.
    Let $l \subset T_2$ be an infinite geodesic with the following properties:
    \begin{itemize}
       \item $O_2 \in l \bigcap l_2$.
       \item $\Gamma_{l} \defeq \Stab_{\Gamma}(l) = \langle e \rangle$
    \end{itemize}
    The existence of such a geodesic is clear by counting considerations:
    Each non identity element $g \in PGL_2(\Q_{q_2})$ fixes
    at most three points on the boundary of $T_2$
    (\cite{Serre:Trees},\ref{sec:BT}(\ref{itm:boundary})),
    but $\Gamma$ is countable and the boundary is not.
    By \ref{sec:BT}(\ref{itm:boundary}) there exists an element $g \in PSL_2(\Z_{q_2})$
    such that $gl_2 = l$.
    The weak approximation theorem \cite{Lubotzky:Book},
    says that $\Sigma \cap PSL_2(\Q_{q_2})$ is dense\footnote{The
          topology here is the $q_2-adic$ topology on $PGL_2(\Q_{q_2})$.
          This is the same as the compact open topology coming from the action
          on the tree $T_2$: two elements are close if their
          restrictions to some large finite subset of $VT_2$ coincide.}
    in $PSL_2(\Q_{q_2})$, and since $PSL_2(\Z_{q_2})$ is open
    $\Gamma \cap PSL_2(\Z_{q_2}) < PSL_2(\Z_{q_2})$ is dense.
    For a given $n$ we can find an element $\gamma_n \in \Gamma$ which is close enough
    to $g$ that $\gamma_nC_n = gC_n \subset l$.
    It follows that $\Gamma_{gC_n} = \gamma_n \Gamma_{C_n} \gamma_n^{-1}$ and that
    the graphs $X_{n} \cong X_{gC_n}$ are isomorphic. In particular, the graphs $X_{gC_n}$
    are also Ramanujan graphs with bounded girth. The intersection
    \begin{equation} \label{eqn:trivial_intersection}
    \bigcap_n \Gamma_{gC_n} = \Gamma_{gC} = \langle e \rangle
    \end{equation}
    is now trivial and therefore the tower of graphs $\{ X_{gC_n} \}_{n \in \N}$
    satisfies all the properties stated in theorem (\ref{thm:main}).
    Note that we have not changed the isomorphism type of the graphs but merely
    the covering maps between them, in order to satisfy equation
    (\ref{eqn:trivial_intersection}).
\qed

\section{Explicit description}
\label{sec:concrete}
    Here we specialize to a very concrete example, which enables us to give an
    explicit description of a tower of Ramanujan graphs in terms of Schreier graphs.

    \begin{theorem} (Compare \cite[theorem 7.4.3]{Lubotzky:Book}).
        Let $q_1,q_2$ be two primes both congruent to $1$ mod $4$, $n$ any integer,
        \begin{eqnarray}
          L(n) & = & \left\{
          \begin{array}{ll}
            PSL_2(\Z/q_2^{n}\Z) & \text{if $q_1$ is a quadratic residue mod $q_2$} \\
            PGL_2(\Z/q_2^{n}\Z) & \text{otherwise} \\
          \end{array} \right. \nonumber \\
        A(n) & \subgroup & L(n) \qquad \text{ the diagonal group. } \nonumber \\
        S(n) & = & \left\{ \left(
        \begin{matrix}
          x_0+x_1\sqrt{-1}  & x_2+x_3\sqrt{-1}  \\
          -x_2+x_3\sqrt{-1} & x_1-x_2\sqrt{-1}
        \end{matrix} \right) \in L(n) \left|
        \begin{array}{l} \sum_{i=1}^4 x_i^2 = q_1 \\
                 0 < x_0 = 1 (\mod 2) \\
                 x_1 = x_2 = x_3 = 0 (\mod 2)
        \end{array} \right. \right \} \nonumber
        \end{eqnarray}
        Then, the set $S$ is symmetric and contains exactly $(q_1+1)$ elements.
        The right Schreier graphs
        \begin{equation} \label{eqn:precise}
             X(L(n),A(n),S(n))
        \end{equation}
        form a tower of $(q_1+1)-$regular Ramanujan graphs with girth $1$
        (i.e. all graphs contain a loop).
    \end{theorem}
    \begin{proof}
    This is just spelling out theorem \ref{thm:main} when $\H = \H_{1,1}$ is the (standard)
    Hamilton quaternion algebra. $\H$ splits at all primes except for $\{2,\infty\}$,
    when $q$ is a prime congruent to $1$ mod $4$ there is a $\sqrt{-1} \in \Q_q$ and
    the splitting is explicitly described by
    \footnote{One can check that this is an isomorphism by solving the linear equations
    for the matrix entries.}:
    \begin{equation} \label{eqn:split}
        \phi_{q}: a+bi+cj+dk \arrow \left( \begin{array}{cc}
                        a+b\sqrt{-1}  & c+d\sqrt{-1} \\
                        -c+d\sqrt{-1} & a-b\sqrt{-1}
                  \end{array} \right)
    \end{equation}
    Let $G = \H^{*}/Z\H^{*}$ be the algebraic group of invertible elements in $\H$
    modulo the center. For the
    lattices \footnote{$\Gamma$ is a uniform lattice in  $PGL_2(\Q_{q_1})$, $\Sigma$
                   is a uniform lattice under its diagonal embedding in
               $PGL_2(\Q_{q_{1}}) \times PGL_2(\Q_{q_{2}})$.
               This means that it acts with a finite number of orbits and
               finite vertex stabilizers on the product $T_1 \times T_2$.}
    $\Gamma$ (resp. $\Sigma$) we take the principal congruence subgroup
    mod $2$, of $G(\Z[1/q_1])$ (resp. $G(\Z[1/q_1,1/q_2])$).
    \begin{eqnarray} \label{eqn:def_Gamma}
        \Gamma & \defeq & \left\{ \left[x_0+x_1i+x_2j+x_3k\right]
            \left| \begin{array}{l}
                    x_i \in \Z                  \\
                \sum_{i=0}^3 x_i^2 = q_1^m {\text{ for some }} m \in \N \\
                x_0 = 1 (\mod 2)                \\
                    x_2=x_3=x_4 = 0 (\mod 2)
            \end{array} \right. \right\}   \nonumber \\
        \Sigma & \defeq & \left\{ \left[x_0+x_1i+x_2j+x_3k\right]
            \left| \begin{array}{l}
                    x_i \in \Z                  \\
                \sum_{i=0}^3 x_i^2 = q_1^mq_2^n
                {\text{ for some }} m,n \in \N  \\
                x_0 = 1 (\mod 2)                \\
                    x_2=x_3=x_4 = 0 (\mod 2)
            \end{array} \right. \right\}
    \end{eqnarray}
    Where the square brackets stand for equivalence class modulo the center $Z\H$.

    The group $\Gamma$ acts freely, transitively and without inversion on the
    tree $T_1$ (see \cite[Lemma 7.4.1]{Lubotzky:Book}),
    this makes all the details of example (\ref{ex:rose}) applicable,
    so we identify the graphs $X_n$ with the Schreier graphs
    \begin{equation} \label{eqn:crude}
       X_n = X (\Gamma/\Gamma(q_2^n),\Gamma_{C_n}/\Gamma(q_2^n),S')
    \end{equation}
    Here $S'$ is the natural set of generators of $\Gamma$, which makes $T_1$ into
    the Cayley graph of $\Gamma$, $\{C_n\}_{n \in \N}$ is an ascending sequence of
    segments of length $2n$ around $O_2$ and $\Gamma(q_2^n) = \Gamma_{B_{T_2}(O_2,n)}$
    is a natural choice for a normal subgroup contained in $\Gamma_{C_n}$.

    The group $\Gamma/\Gamma(q_2^n) = \Gamma/\Gamma_{B_{T_2}(O_2,n)}$ is identified in
    \cite[Remark 7.4.4]{Lubotzky:Book} as the group $L(n)$
    (one can see from \ref{sec:BT}(\ref{itm:spheres}) that it is
    a subgroup of $PGL_2(\Z/q_2^{n}\Z)$). We have seen in the previous section, that
    as long as we are only interested in the isomorphism type of the graphs $X_n$ the
    precise choice of the segments $C_n$ does not matter. We make the choice that will
    give, using \ref{sec:BT}(\ref{itm:spheres}),
    $\Gamma_{C_n}/\Gamma_{B_{T_2}(O_2,n)} = A(n)$.
    To identify $S'$, we choose a base vertex $O_1 \in VT_1$ - the vertex stabilized by
    $PGL_2(\Z_{q_1})$, and let
       \begin{equation}
        S'=\left\{ [x_0+x_1i+x_2j+x_3k] \in \Gamma
        \left|
        \begin{array}{l}
            x_i \in \Z \\
            \sum_{i=1}^4 x_i^2 = q_1 \\
            0 < x_0 = 1 (\mod 2) \\
            x_1 = x_2 = x_3 = 0 (\mod 2)
        \end{array} \right. \right \}
       \end{equation}
    be the symmetric set of generators of $\Gamma$ taking $O_1$ to its
    $(q_1+1)$ neighbors. $S(n)$ is exactly the image of $S'$ in $L(n)$ under
    the map $\psi_{q_1^n} \compos \phi_{q_2}$.

    After making all these identifications equation \ref{eqn:crude} gives
    the desired equation \ref{eqn:precise}.
    From the previous section we know that $\{X_n\}_{n \in \N}$, thus defined
    is a tower of Ramanujan graphs of bounded girth. In order to show that the
    girth is actually $1$ and in order to find the minimal common covering
    of all graphs we must explicitly identify the torus $\langle \gamma,\delta \rangle$.

    Each $q_i$ can be represented as a sum of two squares $q_i = a_i^2 + b_i^2$.
    Without loss of generality we may assume that $a_i$ is positive and odd
    and that $b_i$ is even.
    As our pair of commuting elements we can take $\gamma = a_1+b_1i, \delta = a_2+b_2i$.
    We let $l_i \subset T_i$ be the bi-infinite geodesic which is stabilized (setwise)
    by the diagonal (Cartan) subgroup $A \subset PGL_2(\Q_{q_i})$ then $\gamma$ fixes $l_2$
    pointwise and acts on $l_1$ by translation of length $1$.
    Since by our choice the $C_n$'s are all subsets of $l_2$, the line $l_1$ will be mapped
    into a loop (=circle of length $1$) in each of the graphs $X_n$.

    The minimal common covering of all the $X_n$'s (with the natural covering morphisms)
    is not $T_1$. In order to obtain a tower whose minimal common covering is $T_1$ one has
    to replace $A(n)$ by $g(n) A(n) g(n)^{-1}$ where $g(n) = \psi(q_2^{n}(g))$ is
    the reduction mod $q_2^n$ of any element $g \in PGL_2(\Q_{q_2})$ such that
    $\Gamma_{gl_2} = \langle e \rangle$. For example any element with entries that
    are not algebraic with respect to each other will do, but one can also find
    concrete algebraic examples.
\end{proof}
\noindent
{\bf{Remark: }} Another variant would be to replace $C_n$ by paths of length
            $n$ starting at the vertex $O_2$. The graphs obtained in this way would
        be $X(S(n),B(n),S(n))$ where $B(n) < L(n)$ is the upper triangular
        (Borel) subgroup. These graphs can be identified with the graphs
        coming from the action of $L(n)$ on the projective line
        $\mathbb{P}^{1}(\Z/q_2^{n}\Z)$.
\section {Remarks and open questions}
\begin{itemize}
\item \label{itm:manifolds}
    A similar construction can be carried out for surfaces.
    If we replace $\H$ with a non-definite quaternion algebra which ramifies over $\Q$,
    and $\Sigma$ by a mixed irreducible lattice in a product of a $p-$adic and a
    real Lie group then we can obtain a tower of compact surfaces $S_n$.
    By the Jacket-Langlands correspondence all $S_n$ satisfy Selberg's
    $\lambda_1 \ge 3/16$ (and conjecturally $\lambda_1 \ge 1/4$) theorem,
    and in addition they will all share the same closed geodesic.

    For example pick the quaternion algebra
    $\H = \H_{-2,-3} \defeq \langle 1,i,j,k | ij = -ji = k, i^2 = 2 , j^2=3 \rangle$.
    $\H$ splits at $7$ and at $\infty$ via the splitting
    \begin{equation}
     x_0+x_1i+x2_j+x_3k \arrow
     \left( \begin{matrix} x_0 + x_1\sqrt{2}  &  x_2 + x_3 \sqrt{2} \\
               3 (x_2 -  x_3 \sqrt{2})     &  x_0 - x_1\sqrt{2} \end{matrix} \right)
    \end{equation}
    Let $G$ be the algebraic group $\H^{*} / Z\H^{*}$ defined over $\Q$ and consider the
    lattice $\Lambda \defeq G(\Z[1/7])(2) \in PGL_2(\R) \times PGL_2(\Q_7)$ (i.e. the
    principal congruence subgroup mod $2$ of $G(\Z[1/7])$). Lambda is torsion free.
    The element $\gamma \in \Gamma$ given by
    \begin{equation}
    \gamma \defeq 3 + 2i \sim
           \left( \begin{matrix} 3+2\sqrt{2} & 0 \\ 0 & 3-2\sqrt{2} \end{matrix} \right)
    \end{equation}
    fixes pointwise an infinite geodesic $l_2$ on the tree associated to $PGL_2(\Q_7)$;
    it is also hyperbolic, acting as a translation of length $\log(17 + 12\sqrt{2})$
    on the axis $l_1 = i\R$, as an element of $PGL_2(\R)$ acting on the
    hyperbolic plane $\mathfrak{h}$. We obtain a tower of Riemann surfaces
    $\left\{ \Lambda_{C_n} \bs \mathfrak{h} \right\}_{n}$ where
    $C_n \subset l_2$ are finite segments with length going to infinity.
    All these surfaces will contain a closed geodesic of length $\log(17 + 12\sqrt{2})$.
    As in the combinatorial case the minimal common
    covering of all these surfaces will be $\langle \gamma \rangle \bs \mathfrak{h}$,
    and by twisting
    the covering morphisms we can arrange for the minimal common covering to be
    $\mathfrak{h}$.
\item \label{itm:conj}
    In retrospect the fact that there exist congruence subgroups with non-trivial
    intersection seems obvious. It is interesting to note however that this
    intersection can not be too large. This is the content of the following lemma,
    the proof of which was indicated to me by Alex Lubotzky:
    \begin{lemma} \label{lem:AbelianFG}
        Let $\Gamma < PGL_2(\Q_q)$ be an arithmetic lattice, $\Delta < \Gamma$
        an infinite index torsion free subgroup which is an intersection of
        congruence subgroups.
        Then $\Delta$ is either trivial or $\Z$.
    \end{lemma}
    \begin{proof}
        If $\Delta$ is torsion free then it is free (because it acts freely on
        the tree). If it is not Abelian it must be Zariski dense because $PGL_2$
        does not have any proper, non solvable algebraic subgroups. But $\Delta$
        is closed in the congruence topology on $\Gamma$ so by the strong approximation
        theorem (\cite{Weisfeiler:SAT,Nori:SAT}) $\Delta$ is open in
        the congruence topology on $\Gamma$ and therefore of finite index.
    \end{proof}
    Geometrically this means that for any infinite tower of LPS graphs the
    minimal common covering graph contains at most one circle. In fact, this statement
    will hold for all known constructions of Ramanujan graphs of constant degree. As far
    as I know all these constructions come from congruence subgroups in arithmetic lattices.
    One is lead to ask the following question:
    \begin{question}
        Can the minimal common covering of a tower of Ramanujan graphs
        have more then one circle (i.e. have a non Abelian fundamental group)?
    \end{question}
\end{itemize}
\bibliography{yair}

\end{document}